\documentclass[a4paper,11pt]{amsart}
\usepackage
[
        a4paper,% other options: a3paper, a5paper, etc
        left=2cm,
        right=2cm,
        top=3cm,
        bottom=2.5cm,
]
{geometry}
\usepackage{amsthm,amsmath,amsfonts,amssymb}
\usepackage{dsfont, times}
\usepackage[T1]{fontenc}
\usepackage[utf8]{inputenc}
\usepackage[portuguese]{babel}
\usepackage{hyphenat}
\usepackage{dsfont}           % para exibir melhor o conjuntos dos numeros
\usepackage{mathrsfs}
\usepackage{xcolor}

\usepackage{hyperref}
\hypersetup{
    pdftoolbar=true,
    pdfmenubar=true,
    pdffitwindow=false,
    pdfstartview={FitH},
    pdftitle={Template \LaTeX\ para a Revista Matemática Universitária},
    pdfauthor={Paolo Piccione},
    pdfsubject={},
    pdfkeywords={},
    pdfnewwindow=true,
    colorlinks=true, %set this false for printable version
    linkcolor=blue,
    citecolor=blue,
    urlcolor=blue,
}

\usepackage{graphicx}
\usepackage{amsmath}
\usepackage{amssymb}

\usepackage{fancyhdr}
\def\sloppy{\tolerance=9999 \hfuzz=.5pt \vfuzz=.5pt}
\renewcommand{\contentsline}[4]{\csname nuova#1\endcsname{#2}{#3}{#4}}
\newcommand{\nuovasection}[3]{\medskip\hbox to \hsize{\vbox{\advance\hsize by -1cm\baselineskip=9pt\parfillskip=0pt\leftskip=3.5cm\noindent\hskip -2cm #1\leaders\hbox{.}\hfil\hfil\par}$\,$#2\hfil}}
\newcommand{\nuovasubsection}[3]{\smallskip\hbox to \hsize{\vbox{\advance\hsize by -1cm\baselineskip=10pt\parfillskip=0pt\leftskip=4cm\noindent\hskip -2cm #1\leaders\hbox{.}\hfil\hfil\par}$\,$#2\hfil}}

\def\noi{\noindent}

\newcommand{\demo}{\noindent\textbf{Demonstração.\,}}%começo de uma demonstracao
\newcommand{\fimdemo}{\hspace*{\fill}{$\square$}}  %fim de uma demonstração

\newtheorem{teo}{Teorema}[section]
\newtheorem{defi}{Definição}[section]

\newtheorem{lema}{Lema}[section]

\newtheorem{Obs}{Observação}[section]
\newtheorem{ex}{Exemplo}[section]

\newtheorem{pgt}{Pergunta}[section]

\begin{document}
 %\begin{titlepage}
		
 % ajustar o número da página inicial
 \setcounter{page}{3}

%\hrule\vskip3pt\hrule
\vskip1cm

\title[Alguns teoremas do tipo valor médio]{Alguns teoremas do tipo valor médio: De Lagrange a Malesevic}
 \author[M.Bongarti]{Marcelo Bongarti e German Lozada-Cruz}
 \address{Department of Mathematical Sciences, The University of Memphis, Memphis, TN 38152, USA}
 \email{msbngrti@memphis.edu}
% \thanks{Apoiado pela FAPESP processo: 2013/03866-9}
%\author[G.Lozada-Cruz]{German Lozada-Cruz}
 \address{Departamento de Matemática, Instituto de Biociências, Letras e Ciências Exatas (IBILCE)--Universidade Estadual Paulista (UNESP), São José do Rio Preto, Brasil}
 \email{german.lozada@unesp.br}
%\thanks{Apoiado pela FAPESP processo: 2015/24095-6}
					
\date{07 de janeiro de 2021 (Data da publica\c c\~ao do artigo)}
	 
\maketitle
\thispagestyle{plain}
\begin{small}
%\tableofcontents
\end{small}
				
%\end{titlepage}

\maketitle

\begin{abstract}
Nosso objetivo neste trabalho é apresentar alguns teoremas do tipo valor médio que não são estudados em disciplinas clássicas de cálculo e análise matemática. Trata-se de teoremas simples e de grande aplicabilidade  na análise matemática (por exemplo no estudo de equações funcionais, operadores integrais, etc), matemática computacional, economia e outras áreas.
\end{abstract}

\section{Introdução}

Dos teoremas clássicos  da análise matemática, o \textit{teorema do valor médio} se destaca por sua simplicidade e vasta aplicabilidade. É, sem dúvida, um dos resultados mais conhecidos pela comunidade matemática, e, sem favor algum, um dos tijolos que constituem os alicerces do Cálculo e da Análise Matemática como um todo.

O objetivo deste artigo é apresentar uma gama de outros teoremas do tipo valor médio que somam-se ao clássico e que também são de altíssima aplicabilidade, simplicidade e contribuem para o avanço da Matemática.

Atualmente existem diversas maneiras para abordar o teorema do valor
médio, cada qual ligada a alguma meta específica a qual se pretende
chegar. Neste caso, como queremos um panorama acerca dos diferentes
teoremas do tipo valor médio, utilizaremos a abordagem clássica.

Nossa história começa em 1691 quando Rolle\footnote{Michel Rolle (1652-1719), matemático francês.} usou técnicas do Cálculo diferencial e integral para provar o seguinte resultado, nosso
primeiro teorema do tipo valor médio:

\begin{teo}[Teorema de Rolle]
Seja $f: [a,b] \to \mathbb{R}$ uma função contínua em $[a,b]$ e diferenciável em $(a,b)$. Se $f(a)=f(b)$,
então, existe $c \in (a,b)$ tal que $f'(c)=0$, isto é, a reta
tangente ao gráfico de $f$ no ponto $(c, f(c))$ é horizontal.
\end{teo}

Ainda que Rolle tenha sido merecidamente homenageado pela formalização do resultado, é interessante comentar que Bhaskara II\footnote{Bhaskara Akaria (1114-1185), matemático indiano.}, demonstrou um caso particular do teorema de Rolle muito tempo antes, embora sem pouca ou nenhuma formalidade.

O teorema de Rolle ficou mais conhecido depois que Drobisch\footnote{Moritz Wilhelm Drobisch (1802-1896), matemático alemão.} usou o termo pela primeira vez em 1834, seguido por Bellavitis\footnote{Giusto Bellavitis (1803-1880), matemático italiano} em 1846. Maiores detalhes podem ser encontrados em
\cite{HIST1}.

O teorema do tipo valor médio mais famoso da história é o \emph{teorema do valor médio de Lagrange:}

\begin{teo}[Teorema do valor médio de Lagrange]
Se $f: [a,b] \to \mathbb{R}$ é uma função contínua em $[a, b]$ e derivável em $(a,b)$, então existe $c \in (a,b)$ tal que $$f'(c)=\frac{f(b)-f(a)}{b-a}.$$
\end{teo}

Este resultado foi inicialmente descoberto por Lagrange\footnote{Joseph Louis Lagrange (1736-1813), matemático
italiano.}, que o demonstrou sem, inicialmente, fazer menção ao teorema de Rolle. Entretanto, a dedução mais conhecida tem como ideia principal a aplicação do teorema de Rolle à função auxiliar
\begin{align*}
\varphi(x) = f(x)- \bigg[\frac{f(b)-f(a)}{b-a}(x-a) + f(a)\bigg],
\end{align*}
e esta foi feita por Bonnet\footnote{Pierre Ossian Bonnet (1819-1892),matemático francês.}. Publicamente, o teorema do valor médio de Lagrange foi citado pela primeira vez em um trabalho do renomado
físico Ampére\footnote{André-Marie Ampére (1775-1836), físico francês.}. Para maiores detalhes, a referência \cite{SAHOO} pode ser
consultada.

Geometricamente, o teorema do valor médio de Lagrange diz que existe um ponto $c$ dentro do intervalo $(a, b)$, tal que a reta tangente ao gráfico de $f$ no ponto $(c, f(c))$ é paralela à reta secante que passa pelos pontos $(a, f(a))$ e $(b, f(b))$.

Fisicamente, o teorema do valor médio de Lagrange garante que se uma partícula possui uma trajetória suave $(t, f(t))$ no intervalo de tempo $[a,b]$, então existirá um instante $t_c \in (a,b)$ tal que a velocidade instantânea da partícula (em $t = t_c$) coincide com a velocidade média de todo o percurso.

Motivado por esta aplicação, Cauchy\footnote{Augustine-Louis Cauchy (1789-1857), matemático francês.} se perguntou o que poderia ser dito a respeito de uma partícula de trajetória suave $(f(t),g(t))$
no intervalo de tempo $[a,b]$. Então, Cauchy aplicou o teorema de Rolle à função
$$\varphi(x)=[g(b)-g(a)] f(x)- [f(b)-f(a)]g(x)$$ e o resultado foi o que conhecemos hoje por \emph{Teorema do valor médio de Cauchy}.

\begin{teo}[Teorema do valor médio de Cauchy]
Se $f,g: [a,b] \to \mathbb{R}$ são funções contínuas em $[a, b]$ e deriváveis em $(a,b)$, então existe $c \in (a,b)$ tal que 
\begin{align*}
f'(c)[g(b)-g(a)]=g'(c)[f(b)-f(a)].
\end{align*}
\end{teo}

Geometricamente, o teorema do valor médio de Cauchy garante que dada uma trajetória suave $(f(t),g(t))$ com $t \in [a,b]$, existirá $c \in (a,b)$ de tal forma que a reta tangente à trajetória no ponto $(f(c),g(c))$ é paralela à reta que passa pelos pontos $(f(a), g(a))$ e $(f(b), g(b))$.

\medskip
As demonstrações detalhadas dos teoremas do tipo valor médio discutidos nesta seção podem ser encontradas em \cite[Teorema 6.2.4 e Teorema 6.3.2]{bartle} ou em \cite[Teorema 2.2 e Teorema 2.17]{SAHOO}. O leitor interessado em variações do teorema de Lagrange pode consultar \cite{Lozada-C0}. Para variações e aplicações do teorema de Cauchy, recomendamos \cite{Lozada-C} e \cite{Lozada-C2}

%%%%%%%%%%%%%%%%%%%%%%%%%%%%%
\section{Teorema de Flett e suas variações}\label{def1}

O teorema a seguir é uma versão do clássico teorema do valor médio para integrais e a partir de observações acerca dele que nasce a motivação para o teorema do tipo valor médio que abordaremos a seguir.%do tipo valor médio que iremos abordar neste seção foi motivado por observações acerca do seguinte teorema, que é uma versão do Teorema do valor médio clássico para integrais

\begin{teo}\label{teo:TVMIntegrals}
Se $f:[a,b] \to \mathbb{R}$ é uma função contínua, então existe $\eta \in [a,b]$ tal que $$\int\limits_a^b f(x)dx = f(\eta)(b-a).$$
\end{teo}
A demonstração do Teorema \ref{teo:TVMIntegrals} pode ser encontrada em \cite[Teorema 7.1]{SAHOO}.

\medskip
Considere, agora, uma função $g: [a,b] \to \mathbb{R}$ como no Teorema \ref{teo:TVMIntegrals}, então existe $\xi \in (a,b)$ tal que $$g(\xi)=\dfrac{1}{b-a}\int\limits_a^b g(t)dt.$$

Além disso, se considerarmos uma tal função $g$, contínua, e tal que $$g(a)=0, \ \ \int\limits_a^b g(t)dt=0,$$ e definirmos a função
\begin{align*}
\varphi(x)=\begin{cases}
  \dfrac{1}{x-a}\int\limits_a^x g(t)dt, & x \in (a,b] \\
   0, & x=a,
   \end{cases}
\end{align*}
teremos, então, pelo teorema de Rolle, uma vez que $\varphi$ é contínua em $[a,b]$, derivável em $(a,b)$ e
$\varphi(a)=0=\varphi(b)$, que existe $\xi \in (a,b)$ tal que $\varphi'(\xi)=0$.

Mas, se $x \in (a,b)$, $\varphi'(x)=-\frac{1}{(x-a)^2} \int\limits_a^x g(t)dt + \frac{g(x)}{x-a}$.

Assim, existe $\xi \in (a,b)$ tal que
\begin{equation}
g(\xi)=\frac{1}{\xi-a}\int\limits_a^{\xi}g(t)dt.\label{eq1}
\end{equation}

Note que \eqref{eq1} é exatamente o teorema do valor médio para integrais com a substituição de $b$ por $\xi$. Olhando rapidamente para o teorema fundamental do cálculo, podemos escrever a equação \eqref{eq1} da seguinte forma:

\begin{equation}G'(\xi)=\dfrac{G(\xi)-G(a)}{\xi-a}\label{eq2}\end{equation}
onde $G$ é uma primitiva de $g$, ou seja, $G'=g.$

Neste sentido, dada uma função $g$, é natural nos perguntarmos se podemos trocar a condição de que $\int\limits_a^b g(t)dt =0$ simplesmente por $g(b)=0$. A seguir, como uma consequência do teorema do valor intermediário, vemos que é possível. Estas observações foram feitas por Flett\footnote{Thomas Muirhead Flett (1923-1976), matemático britânico.} e o resultado (de 1958) leva o seu nome como homenagem.

\medskip
O teorema de Flett (veja \cite{FLETT}) é uma variação do teorema de Rolle onde a condição $f(a)=f(b)$ foi substituída por $f'(a)=f'(b)$. Por este motivo, dizemos que o teorema de Flett é um teorema do tipo de
Lagrange com uma condição do tipo Rolle, ou simplesmente, teorema do tipo valor médio com uma condição do tipo Rolle.

\begin{teo}[Teorema do valor médio de Flett \cite{FLETT}]\label{teo:Flett}
Seja $f:[a,b] \to \mathbb{R}$ uma função diferenciável em $[a,b]$ com $f'(a)=f'(b)$. Então, existe $\xi \in (a,b)$ tal que
\begin{align}\label{eqn: TVMFlett}
f'(\xi)=\frac{f(\xi)-f(a)}{\xi-a}.
\end{align}
\end{teo}
\demo Sem perda de generalidade podemos supor $f'(a)=f'(b)=0$, pois caso contrário fazemos
$\psi(x)=f(x)-xf'(a)$ e daí teremos $\psi'(a)=\psi'(b)=0$. Definamos a função $\varphi:[a, b]\to \mathds{R}$ dada por 
\begin{align*}
\varphi(x)=\begin{cases}
 \dfrac{f(x)-f(a)}{x-a}, & x \in (a,b] \\
   f'(a), & x=a.
  \end{cases}
  \end{align*}
A função $\varphi$ é contínua em $[a, b]$, derivável em $(a, b]$ e para $ x \in (a,b)$, $$ \varphi'(x) =
\frac{f'(x)}{x-a}-\frac{\varphi(x)}{x-a}.$$

Observe que, $\varphi(a)=0$. Se $\varphi(b)=0$, pelo teorema de Rolle existe $\xi \in (a,b)$ tal que $\varphi'(\xi)=0$ e o teorema está provado.

\medskip
Suponhamos $\varphi(b)\neq0$. Se $\varphi(b)>0$, segue que 
\begin{align*}
\varphi'(b)=\frac{f'(b)-\varphi(b)}{b-a}=-\frac{\varphi(b)}{b-a} <0.
\end{align*}
Logo, para $\epsilon>0$ suficientemente pequeno existe $x_1 \in (b-\epsilon,b)$ tal que $\varphi(b)<\varphi(x_1)$. Como $\varphi$ é contínua em $(a, x_1)$ e $0=\varphi(a)<\varphi(b)<\varphi(x_1)$, segue do teorema do valor intermediário, que existe $\eta \in (a,x_1)$ tal que $\varphi(\eta)=\varphi(b)$. E então, do teorema de Rolle
aplicado ao intervalo $[\eta, b]$, existe $\xi \in (\eta, b) \subset (a,b)$ tal que $\varphi'(\xi)=0$, isto é,
$f'(\xi)=\frac{f(\xi)-f(a)}{\xi-a}.$

\medskip
O caso $\varphi(b)<0$ é análogo.\fimdemo

\begin{figure}[t]
   \begin{center}
  \includegraphics[scale=0.5]{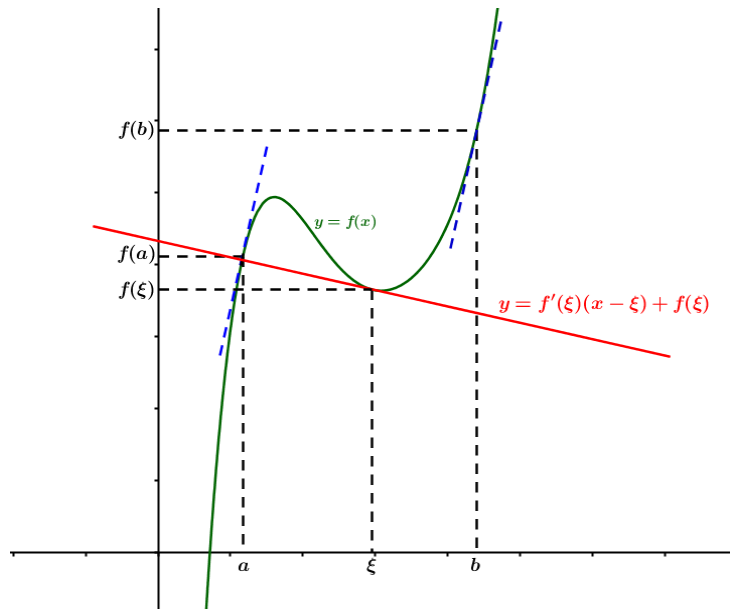}
  \caption{Interpretação geométrica do teorema de Flett}
  \label{figflett}
  \end{center}
\end{figure}

\medskip
Geometricamente, o teorema de Flett diz que se uma curva $(t,f(t))$ é suave no intervalo $[a,b]$ e as retas tangentes nos extremos $(a,f(a))$ e $(b,f(b))$ são paralelas, então, existe um ponto $\xi \in (a,b)$ de modo que a reta tangente ao gráfico de $f$ que passa por $(\xi, f(\xi))$ também passa por $(a,f(a)),$ como podemos ver na Figura \ref{figflett}.

Por outro lado, do ponto de vista cinemático, Flett concluiu que, se as velocidades inicial e final de uma partícula com trajetória $(t,f(t))$ suave no intervalo de tempo $[a,b]$ forem iguais, então, existe um momento $t_j \in (a,b)$ tal que a velocidade instantânea da partícula neste instante, é exatamente a velocidade média do
percurso até o instante $t_j$.

\medskip
Por comodidade, dada uma função $f:[a,b] \to \mathbb{R}$, chamaremos o ponto $\xi\in (a, b)$ tal que satisfaz a conclusão do teorema de Flett simplesmente de \textit{ponto de Flett}.

Um exemplo básico para exemplificar o uso do teorema de Flett pode ser dado quando consideramos a função $f:[-2,2] \to \mathds{R}$ dada por $f(x)=x^3+2x-1$. Como $f$ é um polinômio segue que $f$ é derivável em $[-2,2]$ e usando \eqref{eqn: TVMFlett} vemos fácilmente que um $\xi = 1\in (-2, 2)$ é um ponto de Flett de $f$.

A seguir trataremos brevemente dos resultados apresentados por R. Meyers em 1977. Estes apresentam variações do teorema de Flett. As interpretações geométricas e físicas são análogas às interpretações do teorema de Flett, portanto, não serão verbalmente discutidas nesta seção.

\begin{teo}[{\cite[Teorema $1'$]{MEYER}}] \label{flett2}
Seja $f:[a,b] \longrightarrow \mathbb{R}$ uma função diferenciável em $[a,b]$ com $f'(a)=f'(b)$. Então, existe $\xi \in (a,b)$ tal que
\begin{align*}
f'(\xi)=\frac{f(b)-f(\xi)}{b-\xi}.
\end{align*}
\end{teo}

\begin{figure}[t]
   \begin{center}
  \includegraphics[scale=0.4]{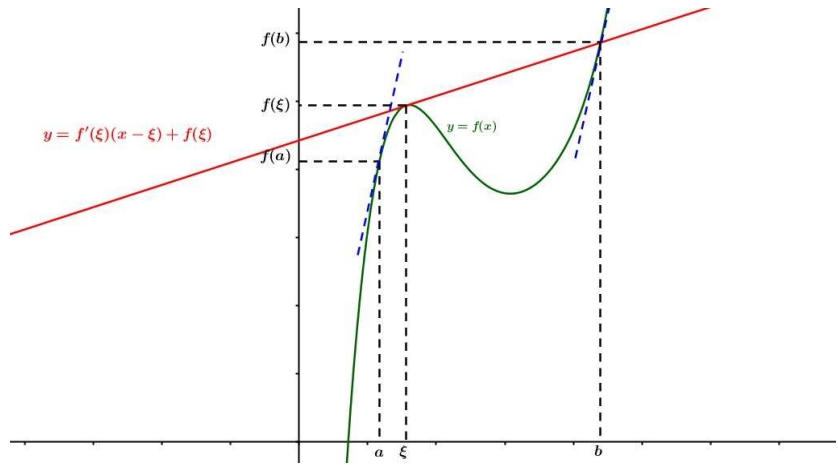}
  \caption{Interpretação geométrica Teorema \ref{flett2}}
  \label{figM1}
  \end{center}
\end{figure}

\begin{teo}[{\cite[Teorema 2]{MEYER}}] 
Seja $f:[a,b] \longrightarrow \mathbb{R}$ uma função diferenciável em $[a,b]$ com $f'(a)=f'(b)$. Então, existe $\xi \in (a,b)$ tal que $$f'(\xi)=\frac{f(b)-f(\xi)}{\xi-a}.$$ \label{flett3}
\end{teo}

\begin{figure}[b]
   \begin{center}
  \includegraphics[scale=0.4]{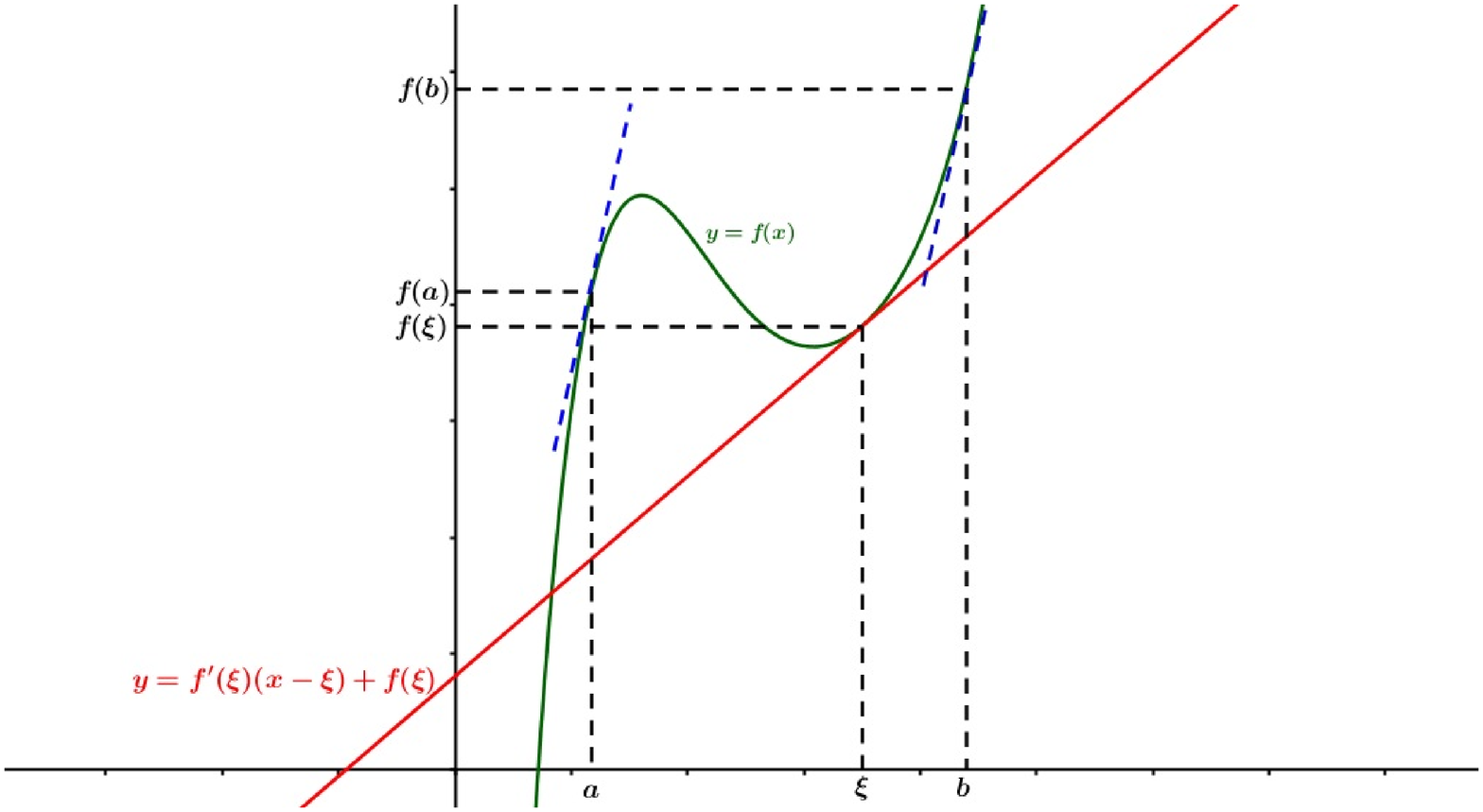}
  \caption{Interpretação geométrica Teorema \ref{flett3}}
  \label{figflett3}
  \end{center}
\end{figure}

\begin{teo}[{\cite[Teorema $2'$]{MEYER}}] Seja $f:[a,b] \longrightarrow \mathbb{R}$ uma função diferenciável em $[a,b]$ com
$f'(a)=f'(b)$. Então, existe $\xi \in (a,b)$ tal que
$$f'(\xi)=\frac{f(\xi)-f(a)}{b-\xi}.$$ \label{flett4}
\end{teo}

\begin{figure}[t]
   \begin{center}
  \includegraphics[scale=0.45]{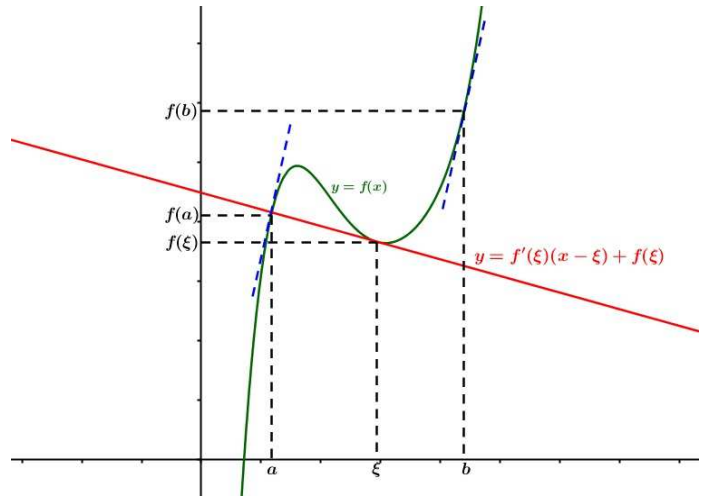}
  \caption{Interpretação geométrica Teorema \ref{flett4}}
  \label{figflett4}
  \end{center}
\end{figure}

\begin{teo}[{\cite[Teorema 3]{MEYER}}] Se $f$ é diferenciável e $f'$ é contínua em $[a,b]$ e
$[f(b)-f(a)][f(b)-f(a)-(b-a)f'(b)]<0$. Então, existe $\xi \in (a,b)$
tal que
$$f'(\xi)=\frac{f(b)-f(a)}{\xi-a}.$$
\end{teo}

\begin{teo}[{\cite[Teorema $3'$]{MEYER}}] Se $f$ é diferenciável e $f'$ é contínua em $[a,b]$ e
$[f(b)-f(a)][f(b)-f(a)-(b-a)f'(a)]<0$. Então, existe $\xi \in (a,b)$
tal que
$$f'(\xi)=\frac{f(b)-f(a)}{b-\xi}.$$
\end{teo}

\begin{teo}[{\cite[Teorema 4]{MEYER}}] Se $f$ é diferenciável e $f'$ é contínua em $[a,b]$ e
$f'(a)[f(b)-f(a)-(b-a)f'(b)]>0$. Então, existe $\xi \in (a,b)$ tal
que
$$f'(\xi)=\frac{f(\xi)-f(a)}{b-a}.$$
\end{teo}

\begin{teo}[{\cite[Teorema $4'$]{MEYER}}] Se $f$ é diferenciável e $f'$ é contínua em $[a,b]$ e
$f'(b)[f(b)-f(a)-(b-a)f'(a)]>0$. Então, existe $\xi \in (a,b)$ tal que
\begin{align*}
f'(\xi)=\frac{f(b)-f(\xi)}{b-a}.
\end{align*}
\end{teo}

%%%%%%%%%%%%%%%%%%%%%%%%%%%%%%
\section{Existência de Pontos de Flett}

O assunto discutido nessa seção tem motivação no exemplo a seguir.

Seja $[a,b]$ um intervalo fechado que contém o $0$ no seu interior e considere a função $f:[a,b] \longrightarrow \mathbb{R}$ dada por $f(x)=|x|$, que não é diferenciável em $x=0$, entretanto, supondo que $a < x < 0$ temos,
\begin{align*}
\dfrac{f(x)-f(a)}{x-a}=\dfrac{|x|-|a|}{x-a}=\dfrac{-x+a}{x-a}
= -1 = f'(x), \  \forall x \in (a,0). 
\end{align*}

Portanto, existem infinitos pontos de Flett em $(a,0) \subset
(a,b)$.

Este exemplo mostra que o conjunto das funções que satisfazem as hipóteses do teorema de Flett está estritamente contido no conjunto das funções que têm um ponto de Flett. Portanto, é natural perguntarmos que outras condições suficientes existem que garantam a existência de pontos de Flett. %Será que existe alguma necessária?

Originalmente, em 1958, T.M. Flett demonstrou que pontos de Flett existem sob as hipóteses de que $f$ seja diferenciável no intervalo fechado $[a,b]$ e que $f'(a)=f'(b)$. Mas esta não é a única.

Os primeiros estudos sobre os resultados de T.M Flett e suas generalizações foram feitos em $1966$ pelo matemático Donald. H. Trahan em 1966 (ver \cite{TRAHAN}). Ele deu uma nova condição para a
existência de um ponto de Flett através de algumas desigualdades, usando uma comparação entre a inclinação da reta secante ao gráfico da função $f:[a,b] \to \mathbb{R}$ passando pelos extremos $(a,f(a))$ e $(b,f(b))$ e a inclinação das retas tangentes ao gráfico passando pelos mesmos.

\medskip
Os seguintes resultados são necessários para a compreensão da condição de Trahan.

\begin{lema}[{\cite[Lema 1]{TRAHAN}}]\label{lema:General-Rolle}
Se $f: [a,b] \to \mathbb{R}$ é uma função contínua, diferenciável em $(a,b]$ e $f'(b)[f(b)-f(a)] \leqslant 0$, então existe $c \in (a,b]$ tal
que $f'(c)=0$.
\end{lema}

\begin{lema}[{\cite[Lema 2]{TRAHAN}}]\label{lema:General-Rolle2}
Se $f: [a,b] \to \mathbb{R}$ é uma função contínua, diferenciável em $(a,b]$ e $f'(b)[f(b)-f(a)] < 0$, existe $c \in (a,b)$ tal que $f'(c)=0$.
\end{lema}

Observe que os Lemas  \ref{lema:General-Rolle} e  \ref{lema:General-Rolle2} são generalizações do Teorema de Rolle.

\begin{teo}[\textbf{Condição de Trahan \cite{TRAHAN}}]\label{teo:Trahan66}
Seja $f:[a,b] \to \mathbb{R}$ uma função diferenciável e tal que
\begin{equation} \label{eqtr}
\big(f'(b)-\tfrac{f(b)-f(a)}{b-a}\big) \big(f'(a)-\tfrac{f(b)-f(a)}{b-a}\big)\geqslant 0.
\end{equation}
Então, existe um ponto de Flett em $(a,b]$.
\end{teo}
\demo Considere a função $\varphi:[a,b] \to \mathbb{R}$ definida por
$$\varphi(x)=\begin{cases}
          \dfrac{f(x)-f(a)}{x-a}, & x \in (a,b] \\
         f'(a) ,& x=a
\end{cases}$$

Observe que $\varphi$ é contínua em $[a,b]$, é diferenciável em
$(a,b]$ e que $\varphi'(b)[\varphi(b)-\varphi(a)] \leqslant 0.$

Logo, pelo Lema \ref{lema:General-Rolle}, existe $\xi \in (a,b)$ tal
que $\varphi'(\xi) =0$, o que significa que
$$f'(\xi)=\frac{f(\xi)-f'(a)}{\xi-a},$$ou seja, $\xi$ é um ponto
de Flett em $(a,b]$. \fimdemo

No exemplo a seguir, temos uma função que não satisfaz a condição de
Flett, mas satisfaz a de Trahan e possui, portanto, um ponto de
Flett.

\begin{ex} \label{trhflt} 
Considere a função $f:[-\tfrac{1}{2},1] \to \mathbb{R}$ dada por $f(x)=x^3$.

Note que $f$ é diferenciável e que  $f'(-\tfrac{1}{2}) \neq f'(1)$,
logo, $f$ não satisfaz a condição de Flett.

No entanto, $f$ satisfaz a condição de Trahan e, portanto, possui um
ponto de Flett, a saber $\xi = \frac{1}{4} \in \left(-\tfrac{1}{2},1\right]$. \label{estricttrahan}\end{ex}

Outra condição suficiente para a existência de um ponto de Flett foi provada por J. Tong em \cite{TONG}. Um ponto interessante desta condição é que Tong só exige a diferenciabilidade de $f$ em $(a,b)$, mas usa os
conceitos de média aritmética de $f$, $\mathscr{M}(f):= \frac{f(a)+ f(b)}{2}$ e média de $f$, $\mathscr{I}(f):=\frac{1}{b-a}\int\limits_a^b f(t) dt$.

\begin{teo}[\textbf{Condição de Tong \cite[Teorema 2]{TONG}}] \label{teot}
Seja $f:[a,b] \to \mathbb{R}$ uma função contínua em
$[a,b]$ e diferenciável em $(a,b)$. Se $\mathscr{M}(f)= \mathscr{I}(f)$ então,
$f$ admite um ponto de Flett em $(a,b)$.
\end{teo}
\demo Basta observar que a função $h$ dada por
$$\begin{array}{ccclc}
   h: &[a,b] & \to & [a,b]\\
     &x & \mapsto & h(x)=\tfrac{f(x)+ f(a)}{2} (x-a)- \int\limits_a^x f(t)dt.
  \end{array}$$
é contínua em $[a, b]$ e diferenciável em $(a, b)$ com derivada $$h'(x)=\frac{1}{2}f'(x)(x-a) + \frac{1}{2} (f(x)+ f(a)) -f(x).$$ Como $h(a)=0$ e $\mathscr{M}(f)= \mathscr{I}(f)$, segue que $h(b)=0$. Agora, usando o teorema de Rolle obtemos a conclusão do eorema. \fimdemo

\medskip
O exemplo a seguir traz uma função que não satisfaz a condição de
Flett, nem de Trahan, mas satisfaz a de Tong.

\begin{ex}Considere a função $f(x)=\arcsin x$
sobre o intervalo $[-1,1].$

Note que $f$ não satisfaz a condição de Flett, nem a de Trahan, pois
não é diferenciável nos extremos. No entanto, um cálculo simples
usando integração por partes garante que $\mathscr{M}(f)=
\mathscr{I}(f),$ o que prova que a função $\arcsin x$ satisfaz a
condição de Tong no intervalo $[-1,1]$ e que, portanto, possui um
ponto de Flett. \end{ex}

\medskip

A terceira condição se deve ao matemático B. Malesevic e é feita em termos de uma função infinitesimal.

Para tanto, seja $f:[a,b] \in \mathbb{R}$ uma função diferenciável em $[a,b]$ e diferenciável um número arbitrário de vezes numa vizinhança à direita do ponto $x=a$.

Considere a expansão de Taylor de ordem um, com resto, dado por
$$f(x)=f(a)+f'(a)(x-a)+\varphi(x)(x-a),$$ onde $\lim\limits_{x \to
a^+} \varphi(x)=0$. Então, definimos a função $\varphi_1: [a,b] \in
\mathbb{R}$ por
\begin{equation} \label{T1}
\varphi_1(x)=
\begin{cases}
\dfrac{f(x)-f(a)}{x-a}-f'(a), & x \in (a,b] \\
0, & x=a.
\end{cases}
\end{equation}

A partir desta função $\varphi_1$ Malesevic demonstrou o seguinte resultado:

\begin{teo}[Condição de Malesevic \cite{MAL}]\label{MAL}
Seja $f:[a,b] \to \mathbb{R}$ é uma função diferenciável e $\varphi_1$ como em \eqref{T1}. Se uma das
seguinte condições
\begin{align*}
{\rm T}_1:&\; \; \varphi_1'(b)\, \varphi_1(b) < 0\;\; e \\
{\rm M}_1:&\; \; \varphi_1'(a)\, \varphi_1(b) < 0
\end{align*}
é satisfeita, então $f$ possui um ponto de Flett. 
\end{teo}
\demo Se a condição ${\rm T}_1$ é satisfeita então $\varphi_1'(b)[\varphi_1(b)-\varphi_1(a)] < 0$. Logo, do Lema \ref{lema:General-Rolle2} existe $\xi_1 \in (a,b)$ tal que $\varphi_1'(\xi_1)=0$, i.e., 
\begin{align*}
\tfrac{1}{\xi_1-a}  \big( f'(\xi_1)- \tfrac{f(\xi_1)-f(a)}{\xi_1-a} \big)=0 \Leftrightarrow f'(\xi_1)= \tfrac{f(\xi_1)-f(a)}{\xi_1-a},
\end{align*}

Agora se a condição ${\rm M}_1$ é satisfeita então $\varphi_1'(a) \big[ \varphi_1 (b)-\varphi_1(a)\big] <0$. Logo, do Corolário 3 do Teorema 2 em \cite{Malesevic} existe $\xi_2\in (a, b)$ tal que $\varphi_1'(\xi_2)=0$, i.e., 
\begin{align*}
\tfrac{1}{\xi_2-a}  \big( f'(\xi_2)- \tfrac{f(\xi_2)-f(a)}{\xi_2-a} \big)=0 \Leftrightarrow f'(\xi_2)= \tfrac{f(\xi_2)-f(a)}{\xi_2-a}.
\end{align*}
\fimdemo

\begin{Obs}
Se ambas as condições ${\rm T}_1$ e ${\rm M}_1$ do {\rm Teorema}  $\ref{MAL}$ são satisfeitas, então existem dois pontos $($distintos$)$ de Flett. $($veja \cite{Malesevic}$)$.
\end{Obs}

\begin{figure}[b]
   \begin{center}
  \includegraphics[scale=0.7]{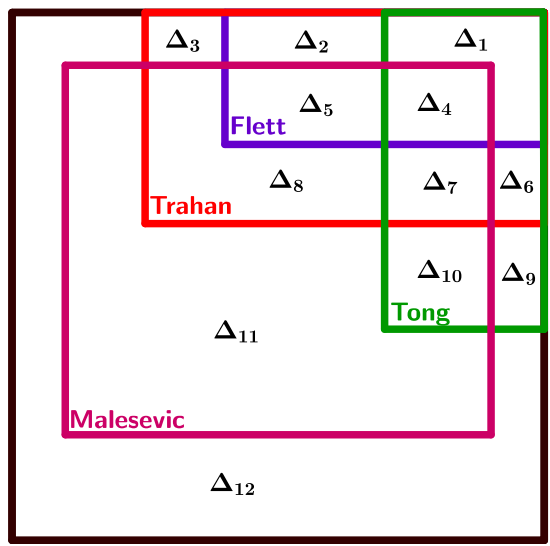}
  \caption{Relação entre as condições de Flett, Tong, Trahan e Malesevic}
  \label{relac}
  \end{center}
\end{figure}

\medskip
A relação entre funções que satisfazem as condições de Flett, Trahan, Tong e Malesevic é representada na Figura \ref{relac}.

Observe que $\Delta_{12} \neq \emptyset$, uma vez que $f(x)= {\rm sgn}(x)$
é uma função que não satisfaz nenhuma das condições, pois não é
diferenciável em $(a,b)$ para qualquer intervalo $[a,b]$ da reta,
que contenha o zero. Entretanto, possui infinitos pontos de Flett.

Analogamente,

\noi {\bf (i)} $f(x)=x^3$, $x \in [-1,1]$ está em $\Delta_1$.

\medskip
\noi {\bf (ii)} $f(x)=\sin(x), x \in \left[-\frac{\pi}{2},\frac{5\pi}{2}\right]$
está em $\Delta_2.$

\medskip
\noi {\bf (iii)} $f(x)=x^3, x \in \left[-\frac{2}{3},1\right]$ está em
$\Delta_3$.

\medskip
\noi {\bf (iv)} $f(x)=\arcsin(x), x \in [-1,1]$ está em $\Delta_{12}$.

\medskip
Prova-se que todos os conjuntos $\Delta_i$, $i=1,2,...,12$ são não vazios, mas a confecção de exemplos para $i=4,5,6,7,8,10$ e $11$ fogem ao objetivo deste trabalho e por isso não serão discutidos aqui.

Maiores detalhes das demonstrações e exemplos vistos nesta seção, sugerimos a referência \cite{MOLNAROVA}.

%%%%%%%%%%%%%%%%%%%%%%%%%%%
\section{Generalizações e Aplicações}

\subsection{Algumas Generalizações}

A seguir trataremos de algumas generalizações e consequências do teorema de Flett. O teorema abaixo foi demonstrado em 1998. Este resultado é uma generalização do teorema de Flett que não exige a condição tipo Rolle, nesse caso o resultado é mais geral. Note que nos dois teoremas a seguir, o caso em que $f'(a)=f'(b)$ é exatamente o teorema de Flett.

\begin{teo}[\textbf{Riedel-Sahoo \cite[Teorema 5.2]{SAHOO}}]
Se $f: [a,b] \to \mathbb{R}$ é diferenciável em $[a,b]$, então, existe $\xi \in (a,b)$ tal que
\begin{align}\label{eqn:Rie-Sa}
f(\xi)-f(a)=(\xi-a)f'(\xi)-\frac{1}{2}\frac{f'(b)-f'(a)}{b-a}(\xi-a)^2.
\end{align}
\end{teo}
\demo Definamos a função $\varphi:[a,b] \to \mathbb{R}$ por
$$\varphi(x)=f(x)-\frac{1}{2}\frac{f'(b)-f'(a)}{b-a}(x-a)^2.$$ Fácilmente vemos que $f$ é diferenciável em $[a, b]$ e $$\varphi'(x)=f'(x)-\frac{f'(b)-f'(a)}{b-a}(x-a).$$ Como $\varphi'(a)=f'(a)=\varphi'(b)$, segue do teorema de Flett que existe $\xi \in (a,b)$ tal que
$$\varphi'(\xi)=\frac{\varphi(\xi)-\varphi(a)}{\xi-a}.$$ Portanto, existe $\xi \in (a,b)$ tal que a equação \eqref{eqn:Rie-Sa}  é satisfeita, o que prova o teorema. \fimdemo

\medskip
Inspirados, então, pela afirmação do Teorema \ref{flett2}, demonstra-se o resultado abaixo.

\begin{teo}[{\cite[Teorema 2.1]{Cakmak}}]
Se $f: [a,b] \to \mathbb{R}$ é diferenciável em $[a,b]$, então, existe $\xi \in (a,b)$ tal que
$$f(b)-f(\xi)=(b-\xi)f'(\xi)+\frac{1}{2}\frac{f'(b)-f'(a)}{b-a}(b-\xi)^2.$$
\end{teo}

O resultado abaixo é também uma generalização do teorema de Flett com outra condição do tipo Rolle $f''(a)=f''(b)$.
\begin{teo}[{\cite[Exercício 5.3.11(b)]{Radulescu}}]\label{teo:f2a=f2b}
Seja $f:[a,b] \to \mathbb{R}$ duas vezes diferenciável e tal que $f''(a)=f''(b)$. Então, existe $\xi \in
(a,b)$ tal que $$f(\xi)-f(a)=(\xi-a)f'(\xi)-\frac{(\xi-a)^2}{2}f''(\xi).$$
\end{teo}

O resultado abaixo é análogo ao Teorema anterior, também com a condição do tipo Rolle $f''(a)=f''(b)$.
\begin{teo}\label{teo:f2da=f2db}
Seja $f:[a,b] \longrightarrow \mathbb{R}$ duas vezes diferenciável e tal que $f''(a)=f''(b)$. Então, existe $\xi \in
(a,b)$ tal que $$f(b)-f(\xi)=(b-\xi)f'(\xi)-\frac{(b-\xi)^2}{2}f''(\xi).$$
\end{teo}

A demonstração do Teorema \ref{teo:f2da=f2db} é análoga à demonstração do Teorema \ref{teo:f2a=f2b} e portanto deixamos como exercício para o leitor. 

Os Teoremas \ref{teo:Flett} e \ref{teo:f2a=f2b} foram generalizados por I. Pawlikowska em \cite{Pawlikowska} para funções $n$ vezes diferenciáveis, com a condição do tipo Rolle $f^{(n)}(a)=f^{(n)}(b)$.

\begin{teo}[{\cite[Lema 2.2]{Pawlikowska}}] 
Seja $f:[a,b] \to \mathbb{R}$ $n$ vezes diferenciável e tal que $f^{(n)}(a)=f^{(n)}(b)$. Então, existe $\xi
\in (a,b)$ tal que $$f(\xi)-f(a)=\sum\limits_{i=1}^{n}\dfrac{(-1)^{i+1}}{i!}(\xi-a)^if^{(i)}(\xi).$$
\end{teo}

%%%%%%%%%%%%%%%%%%%%%%%%%%%%%
\subsection{Algumas Aplicações}
A seguir vamos apresentar algumas aplicações do teorema de Flett. Vamos tratar, principalmente, dos trabalhos feitos por C. Lupu e T. Lupu em \cite{LUPU} e por C. Lupu em \cite{LUPU2}.

\medskip
O objetivo desta seção é apresentar algumas propriedades importantes sobre alguns operadores integrais, como o de Volterra.

Lembremos que $C([0, 1])$ denota o conjunto das funções contínuas reias definidas em $[0, 1]$ e $C^1([0, 1])$ é o conjunto de todas as funções reais continuamente diferenciáveis definidas no mesmo intervalo.

Definamos operadores $T,S:C([0,1])\to C([0,1])$ por 
\begin{align*}
(T\varphi)(t) &=\varphi(t)-\int\limits_{0}^{t}\varphi(x)dx\\
(S\psi)(t) &=t\psi(t)-\int\limits_{0}^{t}x\psi(x)dx.
\end{align*}

As seguintes propriedades valem para os operadores $T$ e $S$.
\begin{teo}[{\cite[Teorema 2.11]{LUPU}}]\label{teo:OpInt}
Se $f,g:[0,1] \to \mathbb{R}$ são funções contínuas, então existem
$\xi_{1},\xi_{2},\xi_{3}\in(0,1)$ tal que
\begin{align*}
\int\limits_{0}^{1}f(x)dx\,(Tg)(\xi_{1}) &=\int\limits_{0}^{1}g(x)dx\, (Tf)(\xi_{1})\\
(Tf)(\xi_{2}) &=(Sf)(\xi_{2})\\
\int\limits_{0}^{1}f(x)dx\,(Sg)(\xi_{3}) &=\int\limits_{0}^{1}g(x)dx\, (Sf)(\xi_{3}).
\end{align*}
\end{teo}

\begin{teo}[{\cite[Teorema 2.12]{LUPU}}]\label{teo:OpInt2}
Se $f,g:[0,1]\to\mathbb{R}$ são funções contínuas, então existem
$\xi_{1},\xi_{2}\in(0,1)$ tal que
\begin{align*}
\int\limits_{0}^{1}(1-x)f(x)dx\, (Tg)(\xi_{1}) &=\int\limits_{0}^{1}(1-x) g(x)\,dx \, (Tf)(\xi_{1})\\
\int\limits_{0}^{1}(1-x)f(x)dx\, (Sg)(\xi_{2}) &=\int\limits_{0}^{1}(1-x) g(x) \,dx\, (Sf)(\xi_{2}).
\end{align*}
\end{teo}
As demonstrações dos Teoremas \ref{teo:OpInt} e \ref{teo:OpInt2} encontram-se com detalhes em \cite[Teorema 2.11 e Teorema 2.12]{LUPU}.

\medskip
O seguinte resultado é a versão equivalente do teorema do valor médio de  Flett para o teorema do valor médio de Cauchy e é criticamente utilizado para obter os resultados seguintes.

\begin{teo}[{\cite[Lema 2.1]{LUPU2}}] \label{flettcauchy}
Sejam $f,g: [a,b] \to \mathbb{R}$ funções diferenciáveis em $[a,b]$ com $g'(x)\neq 0$ para todo $x\in [a, b]$ e
$\frac{f'(a)}{g'(a)}=\frac{f'(b)}{g'(b)}.$
Então, existe $\xi \in (a,b)$ tal que $$\frac{f(\xi)-f(a)}{g(\xi)-g(a)}=\frac{f'(\xi)}{g'(\xi)}.$$
\end{teo}
%\demo Seja $\varphi:[a,b] \to \mathbb{R}$ definida por
%\begin{align*}
%\varphi(x)=\begin{cases}
 %\dfrac{f(x)-f(a)}{g(x)-g(a)}, & x \neq a \\
  %\dfrac{f'(a)}{g'(a)}, & x=a.
  %\end{cases}
%\end{align*}

%Observe que $\varphi$ é diferenciável em $[a,b]$ e, portanto,
%limitada. Se $\varphi$ não atinge seus extremos simultaneamente em
%$a$ e $b$, existe $x_0 \in (a,b)$ onde $\varphi$ atinge seu extremo,
%isto é, onde $\varphi'(x_0)=0$, ou seja,$x_0$ será o número
%procurado.

%Agora, se $\varphi$ atinge seus extremos, respectivamente, em $a$ e
%$b$ temos $\varphi(a) \leqslant \varphi(x) \leqslant \varphi(b)$ ou
%$\varphi(b) \leqslant \varphi(x) \leqslant \varphi(a)$ para todo $x \in
%[a,b].$ 
%
%Suponhamos que o primeiro caso acontece. Note que podemos assumir $g'(x) > 0$ para todo $x \in [a,b]$ uma vez
%que se substituirmos $f$ por $-f$ e $g$ por $-g$ as desigualdades
%acima continuam válidas, donde segue que $f(x) \leqslant
%f(a)+\varphi(b)(g(x)-g(a))$ para todo $x \in [a,b].$

%Assim, $$\dfrac{f(b)-f(x)}{g(b)-g(x)} \geqslant \varphi(b).$$

%Fazendo $x \to b$ temos $\varphi(a) \geqslant \varphi(b)$ e as únicas
%funções que satisfazem $\varphi(a) \leqslant \varphi(x) \leqslant \varphi(b)$
%e $\varphi(a) \geqslant \varphi(b)$ são as funções constantes. Portanto,
%existe $k \in \mathds{R}$ tal que $\varphi \equiv k,$ ou seja,
%$\varphi'\equiv 0,$ o que completa a prova do Teorema.\fimdemo

\medskip
Denotemos por $L^2 ((0,1))$ o espaço vetorial das funções reais quadrado integráveis a Lebesgue sobre $(0,1)$, i.e.,  
\begin{align*}
L^2 ((0,1)) \!=\!\bigg\{\! f: (0,1)\to \mathds{R}\!: f \,\hbox{é Lebesgue mensurável e} {\small \int\limits_a^b f^2(x) dx<\infty}\bigg\}.
\end{align*}
Observemos que para funções contínuas no intervalo $(0, 1)$ a integral de Lebesgue e de Riemann coincidem, logo podemos pensar num primeiro momento no espaço $L^2 ((0,1))$ como sendo o espaço das funções contínuas em $(0, 1)$ cujo quadrado é Riemann integrável em $(0,1)$. Mais ainda, o espaço $L^2 ((0,1))$ é um espaço vetorial normado com a norma $$\|f\|_{L^2((0,1))}=\left(\int_a^b f^2(x) dx\right)^{1/2}, \ \ f \in L^2((0,1)).$$

\begin{defi}[Operador de Volterra]
Sejam $f \in L^2((0,1))$ e  $x \in (0,1)$. Definimos o operador de Volterra, $V$, por
$$\begin{array}{ccclc}
   V: &L^2((0,1)) & \to & L^2((0,1))\\
     &f & \to & V(f)(x)=\int\limits_0^x f(t)dt.
  \end{array}$$
\end{defi}
Sejam $\Psi, \phi:[0,1] \to \mathbb{R}$ funções sendo $\Psi$ contínua e $\phi$ diferenciável com $\phi'(x) \neq 0$ para todo $x \in (0,1)$. Definimos o operador do \textit{tipo Volterra com peso} $$V_{\phi}\Psi(t)=\int\limits_0^t \phi(x)\Psi(x)dx.$$

No que segue, definimos os espaços
\begin{align*}
\mathfrak{C}([a,b])&:=\big\{\phi \in C^1([a,b]), \phi'(x)\neq0, x \in [a,b], \phi(a)\!=\!0\big\}\; \hbox{e}\\
C_{nula}([a,b]) &:= \Big\{ f\in C([a, b]): \int\limits_a^b f(x) dx =0\Big\}.
\end{align*}

\begin{teo}\label{apl0}
Seja $f \in C_{nula}([a,b])$ e $g \in C^1([a,b])$, com $g'(x) \neq 0$ para todo $x \in [a,b]$. Então, existe $\xi \in
(a,b)$ tal que $$V_g f(\xi)=g(a)\cdot Vf(\xi).$$
\end{teo}
\demo Consideremos as funções $\varphi,\eta:[a,b] \to \mathbb{R}$ dadas por 
\begin{align*}
\varphi(t) &=\int\limits_a^t f(x)g(x)dx -g(t)\int\limits_a^t f(x)dx \; \hbox{e}\\
 \eta(t) &=g(t).
\end{align*}
Como $\varphi$ é diferenciável segue que $$\varphi'(t)=f(t)g(t)-\Big(g'(t)\int\limits_a^t f(x)dx +
g(t)f(t)\Big) = -g'(t)\int\limits_a^t f(x)dx.$$ Observe que $\varphi'(a)=0$, assim $\dfrac{\varphi'(a)}{\eta'(a)}=0$.

Por outro lado, $\varphi'(b)=-g'(b)\int\limits_a^b f(x)dx=0$ pois $f \in C_{nula}([a,b]),$ e daí,
$$\dfrac{\varphi'(b)}{\eta'(b)}=0.$$ Logo, do Teorema \ref{flettcauchy} segue que  existe $\xi \in (a,b)$ tal que

$$\frac{\varphi(\xi)-\varphi(a)}{\eta(\xi)-\eta(a)}=\frac{\varphi'(\xi)}{\eta'(\xi)}.$$ Equivalentemente $$\frac{\int\limits_a^{\xi} f(x)g(x)dx - g(\xi)\int\limits_a^{\xi} f(x)dx}{g(\xi)-g(a)}=\frac{-g'(\xi)\int\limits_a^{\xi} f(x)dx}{g'(\xi)}.$$

Disso segue que $$\int\limits_a^{\xi} f(x)g(x)dx = g(a)\int\limits_a^{\xi} f(x)dx.$$ \fimdemo

\begin{teo}\label{apl1}
Se $f,g$ são funções reais contínuas em $[0,1]$ e $\phi \in \mathfrak{C}([0,1])$, então existe $\xi \in (0,1)$ tal que 
\begin{align*}
& V_{\phi}f(\xi)\int\limits_0^1 g(x)dx - V_{\phi} g(\xi)\int\limits_0^1 f(x)dx \\
&\qquad \qquad \qquad \qquad = \phi(0)\bigg(Vf(\xi)\int\limits_0^1 g(x)dx - Vg(\xi)\int\limits_0^1
f(x)dx\bigg).
\end{align*}
\end{teo}
Os detalhes da demonstração do Teorema \ref{apl1} podem ser encontrados em \cite[Teorema 2.4]{LUPU2}. Além disso, algumas observações a respeito deste teorema são pertinentes:

\medskip
\noi ${\bf (i)}$ Se $\phi(0)=0$ então existe $\xi \in (0,1)$ tal que $$V_{\phi}f(\xi)\int\limits_0^1 g(x)dx = V_{\phi} g(\xi)\int\limits_0^1 f(x)dx;$$

Em particular se $\phi(x)=x$, existe $\xi \in(0,1)$ tal que $$\int\limits_0^1 f(x)dx \int\limits_0^{\xi} xg(x)dx
= \int\limits_0^1 g(x)dx \int\limits_0^\xi xf(x)dx.$$

\noi ${\bf (ii)}$ Considere o espaço $L^2$ com peso dado por $$L^2_{\phi}(0,\xi)=\Big\{ u: (0,\xi)\to \mathds{R}:\int_0^\xi u^2(x) \phi(x) \,dx<\infty\Big\}$$ e equipado com a norma $$\|u\|_{L_{\phi}^2(0,\xi)} = \left(\int_0^\xi u^2(x)\phi(x)dx\right)^{1/2}, \ \ u \in  L_{\phi}^2(0,\xi).$$ Substituindo $f$ e $g$ por $f^2$ e $g^2$ temos que existe $\xi \in (0,1)$ tal que
\begin{align*}
& V_{\phi}f^2(\xi)\!\!\int\limits_0^1 g^2(x)dx \!-\! V_{\phi} g^2(\xi)\!\int\limits_0^1 f^2(x)dx \\
&\qquad \qquad \quad =\phi(0) \!\left(\!Vf^2(\xi)\!\!\int\limits_0^1 g^2(x)dx \!-\!\!
Vg^2(\xi)\!\!\int\limits_0^1 f^2(x)dx\!\right),
\end{align*}
ou seja,
\begin{align*}
& ||f||_{L_{\phi}^2(0,\xi)}^2||g||_{L^2(0,1)}^2\!-||g||_{L_{\phi}^2(0,\xi)}^2||f||_{L^2(0,1)}^2\\
&\qquad \qquad \qquad \quad = \phi(0) \big(||f||_{L^2(0,\xi)}^2||g||_{L^2(0,1)}^2 \!-||g||_{L^2(0,\xi)}^2||f||_{L^2(0,1)}^2\big).
\end{align*}
Desta última igualdade, se $\phi(0)=0$, obtemos 
\begin{equation}\label{volt}
||f||_{L_{\phi}^2(0,\xi)}||g||_{L^2(0,1)}=||g||_{L_{\phi}^2(0,\xi)}||f||_{L^2(0,1)}.
\end{equation}
Escrevendo a equação \eqref{volt} da seguinte maneira
\begin{equation}
\frac{||f||_{L_{\phi}^2(0,\xi)}}{||g||_{L_{\phi}^2(0,\xi)}}=\frac{||f||_{L^2(0,1)}}{||g||_{L^2(0,1)}},\label{volt2}
\end{equation}
concluímos a seguinte propriedade interessante: dadas duas funções $f, g$ que tem normas iguais $($ou
proporcionais$)$ em $L^2(0,1)$, e se for dada uma função peso não
constante $\phi$, então existe um número $\xi \in (0,1)$ onde
as normas das funções serão iguais $($ou proporcionais$)$ em $L^2_{\phi}(0,\xi)$.

%%%%%%%%%%%%%%%%%%%%%%%%%
\section{Problemas em aberto}

O estudo de condições necessárias e suficientes para a existência de pontos de Flett (veja seção \ref{def1}) não é, até onde sabemos, completo. Recorde que a função $f(x) = \mbox{sgn}(x)$ pertence ao conjunto $\Delta_{12}$ (veja a Figura \ref{relac}), portanto, não satisfaz nenhuma das condições discutidas anteriormente. No entanto, possui infinitos pontos de Flett. Esta observação, sozinha, torna natural três perguntas, ainda não respondidas na literatura:

\begin{pgt}
Além das apresentadas neste trabalho, existem outras condições suficientes para a existência de pontos de Flett?
\end{pgt}

\begin{pgt}
Existe uma condição necessária para a existência de pontos de Flett?
\end{pgt}

\begin{pgt}
Assumindo que uma função possua pelo menos um ponto de Flett. Sob quais condições este seria único?
\end{pgt}

\medskip
\noi {\bf Agradecimentos:} 
Os autores agradecem ao parecerista pelos comentários e observações, que ajudaram a melhorar de maneira significativa a apresentação deste trabalho. 

Este trabalho é fruto da Iniciação Científica realizada pelo primeiro autor, quem agradece o apoio da FAPESP através do Processo 2013/03866-9.

%%%%%%%%%%%%%%%%%%%%%%%%%%%%%%%

\end{document}